\documentclass[lang=american]{ems-mag-author}

\usepackage{tikz}
\usetikzlibrary{arrows.meta,positioning,fit,backgrounds}
\usepackage[most]{tcolorbox}

\hypersetup{
  pdftitle={Research profiles and the mathematical landscape},
  pdfauthor={Helmut H. W. Hofer}
}

\newcommand{\Profile}{\mathcal P}
\newcommand{\Landscape}{\mathcal L}
\newcommand{\Model}{\mathcal M}

\newcommand{\Probe}{\operatorname{Probe}}

\newtcolorbox{profilebox}{
  breakable,
  fontupper=\small,
  colback=blue!4,
  colframe=blue!35,
  boxrule=.4pt,
  arc=0pt,
  left=5pt,
  right=5pt,
  top=5pt,
  bottom=5pt,
  before skip=7pt,
  after skip=7pt
}

\begin{document}

% Preprint page style: retain the Draft17 layout, with plain page numbers.
\makeatletter
\renewcommand*\ps@plain{%
  \let\@evenhead\@empty
  \let\@oddhead\@empty
  \def\@evenfoot{\hfil{\footnotesize\thepage}\hfil}%
  \let\@oddfoot\@evenfoot
}
\let\ps@doipage\ps@plain
\makeatother
\pagestyle{plain}

\title{Research profiles and the mathematical landscape}
\subtitle{A grassroots proposal for mathematics with AI}
\author{Helmut H. W. Hofer}
\maketitle

\begin{abstract}
Repeated mathematical work with AI leaves behind corrections, failed
arguments, and sometimes useful ways of asking questions. I propose keeping
these in a research profile that the mathematician develops and controls.
The profile marks parts of the mathematical landscape according to the
researcher's understanding and makes those markings available to later work
with AI. This suggests experiments with individual profiles, shared readings
of papers, and graduate education. It also raises mathematical questions
about the landscape itself: its connections, distances, and changes under
discovery. I describe two examples, searching for branches of mathematics
with a common base and recovering useful questions from the history of
Floer theory. The aim is to give mathematicians something to try and a part
in deciding how mathematical AI develops.
\end{abstract}

\section{The profile comes first}

Much of the discussion of AI in mathematics leaves mathematicians waiting
to discover how much of their work the machines will take over. I would
prefer that we take part in deciding what these machines become. There is
an immediate way to begin, in our own mathematical work.

I tried AI in my mathematical work in 2025 and was disappointed. It helped
with routine tasks, but I could not use it as a serious research collaborator.
That changed with a newer model. In late May 2026, I began a sustained
collaboration on a project connecting different areas of mathematics. I had
already spent several years developing the mathematical direction. Working
with AI led me to keep a record of corrections, failed approaches, and
questions that helped the collaboration. This became what I call a research
profile \cite{Hofer2026}.

In my work, substantial profile entries about a research problem and its
surrounding area have made a difference. They bring into the collaboration
my understanding of the problem, possible connections, failed approaches,
and reasons for pursuing one direction rather than another. Later
conversations can build on earlier work, and I have to make some corrections
less often.

The methodological part records how I want to work and what experience has
taught us about doing so. An entry might arise because the model repeatedly replaces
a difficult proposed connection by an easier neighboring problem. The entry
would preserve the correction, the reason for it, and a question to ask
when the same difficulty appears again. We can then see whether it helps.
An instruction that ceases to help should be changed or removed.

I regard a developed profile as money in the bank. It contains understanding
accumulated through work and may improve what I can do in collaboration
with AI. This raises the question whether the privacy promised by AI
companies protects the profile. We should keep that question in mind.

Publication releases our work at distinct moments, and that work can add
to the capabilities of later AI models. If profiles become widely used
and information from them finds its way into future models, by agreement
or through leakage, many boats might continuously lift the tide. But that
is a different story.

A useful picture is that the profile marks parts of the mathematical
landscape according to the researcher's developing understanding. Certain
constructions acquire importance, an obstruction rules out an attractive
route, or a known theorem takes on a new role. These markings express
individual judgment. Two mathematicians reading the same paper may mark
the landscape differently. With their respective profiles available, the
same model answering the same question may attend to different structures
and suggest different approaches. Different answers need not be better
answers. Whether the markings help us do mathematics is a question we can
investigate.

The profile is connected to the mathematical record of each project. We keep
definitions, calculations, manuscript versions, and open proof obligations
in that record. The profile preserves the understanding and ways of working
we want to carry into other problems, with references back to the arguments
that support them. This allows it to remain concise without losing the
mathematics from which it grew.

I think of such a profile as an \emph{adaptive probe}. It affects what we ask
of mathematics and where we look next; experience changes the probe. Its
reach need not be confined to the researcher's present field. A question
arising in one area may make something in a distant area newly visible.
This is where the individual profile begins to meet the larger mathematical
landscape.
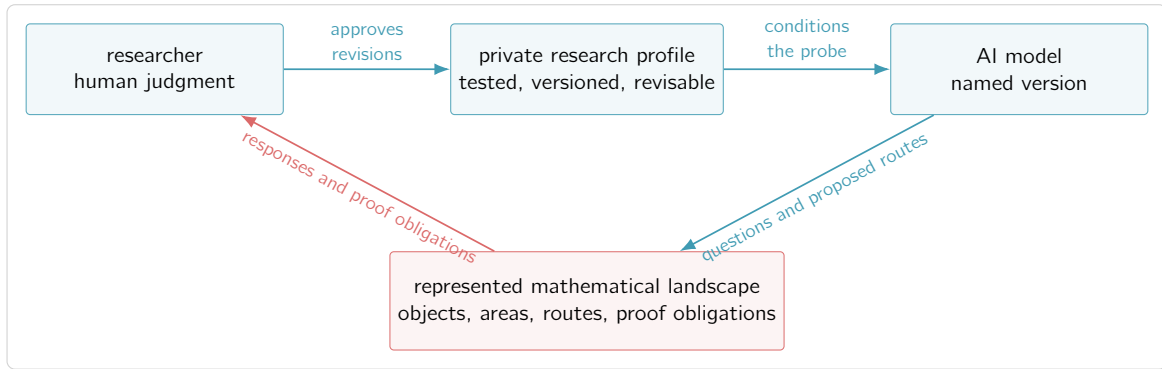
\begin{figure*}[!t]
\centering
\begin{tikzpicture}[
  >=Latex,
  node distance=18mm and 22mm,
  object/.style={draw=blue!65, fill=blue!5, rounded corners=1.5pt,
    minimum width=34mm, minimum height=12mm, align=center, inner sep=3pt},
  landscape/.style={draw=red!65, fill=red!5, rounded corners=1.5pt,
    minimum width=48mm, minimum height=13mm, align=center, inner sep=3pt},
  flow/.style={->, line width=.65pt, blue!75},
  return/.style={->, line width=.65pt, red!70},
  every node/.style={font=\small}
]
  \node[object] (human) {researcher\\human judgment};
  \node[object, right=of human] (profile) {private research profile\\tested, versioned, revisable};
  \node[object, right=of profile] (model) {AI model\\named version};
  \node[landscape, below=18mm of profile] (landscape) {represented mathematical landscape\\objects, areas, routes, proof obligations};

  \draw[flow] (human) -- node[above, align=center, font=\scriptsize] {approves\\revisions} (profile);
  \draw[flow] (profile) -- node[above, align=center, font=\scriptsize] {conditions\\the probe} (model);
  \draw[flow] (model) -- node[sloped, below, font=\scriptsize] {questions and proposed routes} (landscape);
  \draw[return] (landscape) -- node[sloped, below, font=\scriptsize] {responses and proof obligations} (human);

  \begin{scope}[on background layer]
    \node[draw=black!18, rounded corners=3pt, fit=(human)(profile)(model)(landscape),
      inner sep=7pt] {};
  \end{scope}
\end{tikzpicture}
\caption{Work with the model produces proposed routes and new questions. The
researcher judges what to preserve in the profile. Mathematical results
enter the shared landscape through checking, explanation, and communication.}
\label{fig:loop}
\end{figure*}

\section{Related work}

There are precedents for several parts of this proposal. Weisberg and
Muldoon study the division of cognitive labor among researchers moving
over epistemic landscapes \cite{WeisbergMuldoon2009}. Bibliographic data
have been used to represent a space of knowledge and researchers' movement
through it \cite{Singh2023}. Human--AI coevolution concerns the ways people
and algorithmic systems change one another \cite{Pedreschi2024}.

Recent systems also address the problem of retaining scientific experience.
DivAlign uses public research histories and gaps to propose diverse
directions \cite{Xu2026}. Sibyl-AutoResearch and EvoScientist preserve
successful and unsuccessful trials for later work \cite{Wang2026,Lyu2026}.
In materials science, Liu et al. develop a portable memory containing facts,
skills, protocols, and checks \cite{Liu2026}. These are relevant attempts
to make earlier work useful again.

My interest is in a profile developed by a mathematician through sustained
work, under that mathematician's control. It brings an individual
understanding of mathematics into the collaboration. Comparing such
profiles may help us understand why a connection is visible from one
viewpoint and difficult to see from another. The further possibility is
that these differences, and the landscape they reveal, admit a useful
mathematical description.

\section{From profiles to a landscape}

AI can supply considerable muscle within a developed area of mathematics.
It can search for a result, compare conventions, carry out calculations,
write code, or attack a proposed proof. But research also takes us to the
boundary of an area, where familiar methods no longer suffice, and sometimes
to a connection between areas whose relation was not apparent before.

It is useful to distinguish these scales. Locally, we work within a dense
body of examples, methods, and theorems. At an intermediate scale, we may
reorganize a subject or move a method into a neighboring theory. At a global
scale, we try to build bridges between whole areas. We can also ask about
the landscape formed by those areas and bridges, and how discovery changes
it. A profile should help us move between these scales without losing the
question that made the movement worthwhile.

The contribution of the methodological part is less clear. In our small
informal comparisons, the model often explored a specified mathematical
neighborhood about as well without these instructions, given the same
information and prompt. They may matter more when we have to find relations
between areas or keep a larger research direction in view. This needs to
be tested.

For a proposed bridge, we need to identify the objects on both sides,
what corresponds to what, and which structures survive the passage. The
hypotheses and the point where the construction might fail matter as much
as the initial resemblance. A first bridge may be sparse. Examples,
converse constructions, and intermediate results can gradually fill the
region around it until a new area begins to form:
\[
 \text{bridge}\longrightarrow\text{new direction}
 \longrightarrow\text{denser connections}\longrightarrow\text{new area}.
\]

How might we describe this mathematically? I have called it a
\emph{category-type structure}. Objects and relations occur at several
scales. A theorem may be used in a proof, specialize another theorem, or
obstruct a proposed construction. A group of results may eventually be
treated as one object when we look at a larger scale. Often several
ingredients must act together, so a graph with only pairwise arrows may
lose something essential. Category theory and hypergraphs offer possible
languages; neither has yet been shown to give the description we need.

The relations would also carry different weights. A result may be reliable
but irrelevant to the current question, or provisional but strategically
important. The description must allow these distinctions and their changes
over time. I would let examples and experiments determine how much
mathematical structure is warranted.

\subsection{Distances that depend on the researcher and the model}

A theorem can be close for one mathematician and remote for another.
Learning a technique, recognizing a connection, and proving an implication
are different ways of crossing a distance. The passage may be easier in one
direction than in the other. Models, too, find different routes, and a
profile can change which routes a model finds. This suggests quantities
such as
\[
 d_{\mathrm{Gromov}}(A,B),\quad d_{\mathrm{Human}}(A,B),\quad d_{\Model}(A,B),\quad
 d_{\Model,\Profile}(A,B).
\]
The first is imagined relative to Gromov's understanding, the second to
another researcher's. These symbols express questions for a future theory.
We have not defined metrics or established their axioms. We might need several quantities,
depending on the task, the resources available, and the state of knowledge.
The differences between researchers and models are part of what we would
like to understand.

\section{The profile changes through use}

Let $\Landscape_t$ denote the represented landscape at time $t$,
$\Profile_{i,t}$ the profile of researcher $i$, and $\Model_v$ a particular
model version. Schematically, asking a question gives
\[
 (y_t,\Gamma_t,\Omega_t)
 =\Probe_{\Profile_{i,t},\Model_v}(\Landscape_t,q_t).
\]
Here $y_t$ is a proposed answer, $\Gamma_t$ records possible routes and
connections, and $\Omega_t$ contains tests and unresolved proof obligations.
The profile helps determine how the question is pursued. After examining
the work, the researcher may change the profile:
\[
 \Profile_{i,t+1}=U_i(\Profile_{i,t},q_t,y_t,
 \text{human judgment}).
\]
The new profile affects subsequent work. A checked result may also change
the public landscape and, through it, other researchers' profiles. Figure
\ref{fig:loop} shows this continuing exchange.

Several kinds of change must be kept apart. A conversation changes the
immediate context. Continued work changes the researcher and the profile.
Checking and communicating results changes shared mathematics. The model
itself changes through separate training and development; keeping an external
profile does not mean that its parameters change after each conversation.

Any description of the process must retain enough history to account for
these effects. The landscape alone will generally not determine the next
move. A probabilistic model might ask how long it takes to reach a new area,
whether a route survives a change of model, or how an erroneous profile
entry affects later work. An answer could tell us which features of the history matter for
subsequent work.

\section{Latent highways}

An AI model can make a connection look unexpectedly easy. Its training,
architecture, and tools may have placed close together things that are far
apart in a mathematician's working knowledge. I call such a route a
\emph{model-relative mathematical highway}. Using the tentative distance
notation, we might picture
\[
 d_{\mathrm{Human}}(A,B)\gg 1,\quad
 d_{\mathrm{Gromov}}(A,B)<1,\quad
 d_{\Model}(A,B)\ll 1.
\]
A surprising answer may follow a route already available to the model.
A profile may help disclose another route by supplying the researcher's
orientation. In a further case, the work constructs new intermediate
mathematics and establishes a bridge. These possibilities can interact;
the final theorem alone does not tell us how it was found.

Polyfold theory gives me a personal example. Compactified moduli spaces of
pseudoholomorphic curves have enough geometric structure to support
invariants, although breaking, gluing, and failures of transversality put
them outside the usual manifold and Fredholm frameworks. Polyfold and
Fredholm theory provide a setting in which these problems can be perturbed
and their solution spaces studied by algebraic and topological methods
\cite{HWZ2021}.

I had often wondered whether M-polyfolds, beautiful geometric objects in
their own right, might be useful elsewhere. It seemed only ``fair'' that
they should. While writing this article, I asked, with my profile available,
whether a calculus on spaces of varying dimension could be useful in AI.
The answer brought to my attention a paper on PolyNODEs that I had not known
\cite{AhagEtAl2026}. The authors pass from neural ODEs on Euclidean spaces
to manifolds and then to M-polyfolds, where a state can move continuously
through dimensional bottlenecks. They construct examples and give
proof-of-concept experiments in reconstruction and classification.

The application already existed. The useful event was that a question
arising from my own mathematics found its way to it. This may be an example
of a highway available to the model. I can now ask a more pointed question:
which of the structures that mattered in developing polyfold theory might
be useful in AI? The profile can retain the reasons for asking, including
insights that would take time to reconstruct from the published theory.

The construction of a complex structure on $S^6$ gives another example of
what remains to be done after an AI announcement. Alp\"oge's initial
account and its long AI-generated development were followed by Alexeev's
Lean repository and Feng's explanatory workflow
\cite{Alpoge2026,Alexeev2026,Feng2026}. On 29 September 2026, Engel posted
a self-contained geometric construction and proof of a family of compact
complex threefolds diffeomorphic to $S^6$ \cite{Engel2026}. He reports using
ChatGPT in exploration and drafting, checking the mathematics independently,
and taking responsibility for it. His reorganization helps explain the
construction and suggests further questions. Finding the result and making
its ideas available for further mathematics are both substantial work.

OpenAI's release of mathematical manuscripts on 6 October 2026 makes this
task much larger \cite{OpenAIMath2026}. On 22 October, the catalogue lists
719 manuscripts in 372 families and reports Lean formalizations for about
42\% of its top-line results. The repository distinguishes stages of
verification. Its 7 October history records three withdrawals and revisions
to fourteen other manuscripts. We should examine the arguments and the correspondence
between formal statements and the claims made for them, while asking what
the proofs teach us. Section~\ref{sec:profilework} proposes doing this in
groups, with graduate students participating. The work could contribute
both to our understanding of the papers and to the beginnings of individual
research profiles.

\section{Profiles in a mathematical community}

The collaboration involves a researcher, a profile, and a model:
\[
 (\text{researcher},\ \text{profile},\ \text{AI model}).
\]
Continued work changes the first two, while new versions of the third bring
different strengths and weaknesses. With many researchers, we can ask which
parts of their experience are useful to others and what happens when their
different views meet.

Much of a profile will remain private. In my own project, it would reveal
which areas I am trying to connect, why I think the connection matters,
and where the serious difficulties lie. That is mathematical information
before there is a theorem. There is also a simpler consideration:
mathematics is fun, and I enjoy finding the bridge. I want help with that
work while retaining the freedom to pursue it myself.

The researcher should control the files and decide which parts to share.
A shared component might contain a useful checking procedure, a failed
construction with its cause, or questions that brought an overlooked
structure into view. Others could try it without receiving the private
history behind it. We still need to establish what protection particular
systems offer for unpublished work, including what their outputs can
reveal. The profile proposal does not by itself solve that problem.

We could begin comparisons with a capable model given an ordinary prompt,
then supply generic guidance, an individual profile, or selected components
from several profiles. Keep the mathematical sources, model, and resources
comparable. Ask whether the profile changes the questions and whether the
change is useful: does it expose an obstruction, recover a dormant idea,
or reduce the work needed to check an answer? Retain the unsuccessful
cases as well. Otherwise repeated attempts can produce an attractive
example with little evidence that the profile helped.

These experiments need not require pooling private files. Participants can
agree on questions and release selected answers. Dated profile copies make
later comparisons possible. One should also look for harm: an instruction
may preserve an error or repeatedly send different people toward the same
fashionable ideas. A useful collection of profiles should make room for
different mathematics and for disagreement about where to go next.

\section{Human direction}

Mathematicians should have a part in choosing the questions put to these
machines. Computation, energy, and the time needed to check and understand
an argument are finite. We need judgments about what is worth doing, and
people able to explain and contest those judgments.

Questions from physics, biology, engineering, or medicine can reveal a
missing connection or force us to develop new mathematics. AI may help us
learn enough of another subject to recognize such a question and begin
working on it. Conversely, a mathematical structure may change how a
scientist formulates a problem. A richer view of these dependencies could
help us see opportunities that would otherwise remain hidden.

This leaves ample room for pure mathematics whose uses we cannot foresee.
No measure of immediate usefulness should decide its future. Different
mathematicians should be free to pursue different judgments of importance.
We also need colleagues with enough mathematical and technical knowledge
to examine AI systems, question their claims, and influence their
development. Using AI to help with a check does not relieve us of that
responsibility.

\section{A grassroots program}

We need not wait for a platform or an agreed theory of profiles. Researchers
can begin with small profiles in their own work, compare experiences, and
share the parts they find useful. Groups can investigate particular
questions without disclosing their private records. This will give us
examples and failures from which a more systematic study can grow.

The mathematical description should develop alongside these experiments.
If an apparent distance or pattern of interaction keeps recurring, we can
ask what definition would capture it and what that definition would allow
us to prove. We may have to change the language several times. This seems
preferable to building an elaborate theory around objects nobody has used.

The immediate invitation is to take one piece of mathematical work, begin
a profile, and see what happens. The following instruction offers a start.

\section{An instruction to begin}

A reader can give this instruction to an AI and change it as experience
suggests. One can also give this whole article to an AI and begin a
dialogue about starting one's own profile.

\begin{profilebox}
\textbf{Beginning a private research profile}

\smallskip
Help me develop a research profile through our mathematical work. Begin
with a problem, argument, or paper I supply. Preserve the question and its
intended direction, especially when a nearby problem would be easier.

Keep a record of substantive corrections and failed approaches, with the
reasons for them. When we propose a connection, identify the objects,
what corresponds to what, the hypotheses, and the likely obstruction.
Distinguish a theorem from a conjecture, an analogy, or a calculation still
requiring proof. If I ask whether you agree, check the argument independently.

Suggest profile entries when an episode teaches us something useful for
later work. Show me the episode, your interpretation, the proposed entry,
and a way to test it. Say where it applies and when it should be removed.
Distinguish general guidance, my own ways of working, and a correction
needed only for this model or project. Add an entry only with my approval.

Keep the detailed mathematics in the project record and connect the profile
to it. My unpublished work and full profile remain private. Prepare a
shared component only when I request it, and show me what it contains.
\end{profilebox}

\section{Profiles through mathematical work}
\label{sec:profilework}

A paper worth understanding offers a place to begin. Read it with an AI,
work through examples, and ask what makes the argument work. One reader may
see an analytic device, another a geometric construction, and a third an
organization of familiar ingredients that suggests a different use. Each
develops a view of the paper. The profile can preserve the questions and
insights through which that view emerged, and make them available when
another problem brings the paper back into consideration.

Understanding the organization may suggest an abstraction. Which
ingredients can be replaced? Which relations between them must survive?
A failed extension may reveal a dependency that the first reading missed.
Sometimes the right abstract setting is already available. Perhaps, in
other cases, describing the organization will require mathematics that
does not yet exist. That possibility interests me. It should be pursued
through examples where a new description could expose a relation or lead
to a construction we could not otherwise see.

The OpenAI release gives groups of mathematicians an opportunity to develop
research profiles through shared work. As they check and discuss a paper,
each participant can develop a view of its mathematics and begin a profile
of their own. At the same time, the group can experiment with ``profiles
over profiles'' for collaboration. Such a profile could help bring the
individual profiles into play together: which perspectives are useful for
the present question, how they complement or challenge one another, and
what the group learns from their interaction. Participants could then carry
selected parts of that shared experience into their own profiles. The
individual profiles and the profile of the collaboration would develop
together.

This seems to me particularly useful for graduate students. They can learn
the mathematics while beginning a profile with entries that have content
and a history. A student might discover that a simple example exposes a
hidden assumption, that an apparently essential hypothesis is used at only
one step, or that two proofs organize the same ingredients differently.
The detailed argument belongs in the paper's study record. The student's
profile can retain the resulting understanding and useful questions, with
links back to that record. An unresolved point can be an entry too, provided
its status remains clear. There is no need to wait for a complete account
of one's mathematical taste before starting.

The participants should first have room to develop their own views. A
focused shared component could then collect those views, the grounds for
them, and the disagreements. Each reader could add it to an individual
profile and make use of it in a different way. Comparing views may reveal
something none of the readers initially saw. We should try this and find
out whether the shared component helps on later problems, beyond what the
same sources and ordinary AI assistance would provide.

An advisor could also give a PhD student selected parts of a personal
profile: a way of testing a proposed connection, questions that have proved
useful, or a mathematical episode that changed the advisor's view. The
student can try these, revise them, and gradually develop a profile of
their own. The important educational question is what the student comes
to understand and can subsequently do independently.

A researcher acquires this AI muscle over time. A strong profile gives
the collaboration access to accumulated understanding and correction.
Its parts need not all be active at once. The researcher may call one up;
the AI may suggest that a relevant part is available, or activate it under
an arrangement the researcher has chosen. Choosing the right occasion
matters. A good question asked at the wrong time can send us away from
the work we need to do. Activation, reminders, and the revision of borrowed
components are possibilities to test and refine.

The fork and jungle examples below give two more ways to develop a profile.
They ask what we can learn from the branching of theories and from the
difficult circumstances in which new mathematics was found.

\section{Forks and the mathematical jungle}

Here are two instructions readers can try. The first asks for mathematical
developments that have a common base and later separate. The second uses
the history of a discovery to recover questions worth asking in a present
problem. We have called them \emph{forky prompting} and \emph{historical
prompting}.

There are precedents for both: literature-based discovery, structure
mapping, and analogical prompting for the first
\cite{Swanson1986,Gentner1983,Yasunaga2024}; rational reconstruction and
the study of discovery through historical cases for the second
\cite{Lakatos1976,GrasshoffMay1995}. The instructions below make particular
uses of these ideas within a research profile.

The motivation for this instruction predates AI. In 2002, Lai-Sang Young,
Dennis Sullivan, Misha Lyubich, Dusa McDuff, and I submitted an NSF Focused
Research Group proposal, \emph{Symplectic Geometry and Qualitative Dynamics
in Low Dimensions} \cite{FRG2002}. Both subjects had grown from the study
of classical Hamiltonian motion, but had developed powerful methods with
comparatively little interaction. We proposed bringing them together and
explicitly anticipated the ``creation of a new field `Symplectic Dynamics'
with highly integrated ideas and methods from both disciplines.''

The proposal asked us to identify counterparts of core dynamical ideas in
symplectic geometry and conversely. There were already concrete openings:
where holomorphic curves provide surfaces of section for a Reeb flow,
return maps allow methods of surface dynamics to enter. The larger task
was to develop these connections and find others. The NSF did not fund
the proposal, considering it too speculative. To be fair, the opinions of
the referees were as polarized as today's politics. Symplectic dynamics has
since become an established field. Its achievements include the Conley
conjecture for closed symplectically aspherical manifolds
\cite{Ginzburg2010} and smooth closing lemmas in low-dimensional
Hamiltonian dynamics \cite{Irie2025}. The direction that motivated our
proposal is the one behind the forky instruction: return to a common base,
examine what developed on its different branches, and look for further
bridges. AI gives us another way to pursue it.

A \emph{fork object} is a mathematical object,
construction, equation, or problem that is a common starting point for
developments which later separate into different languages, techniques, and
research cultures.  Schematically,
\[
\begin{array}{ccccc}
&&X&&\\[-2pt]
&\swarrow&&\searrow&\\[-2pt]
A_1&\longrightarrow&A_2\longrightarrow\cdots
&\qquad B_1\longrightarrow B_2\longrightarrow\cdots .
\end{array}
\]
At the starting level the two branches share $X$, so there is the trivial
bridge $X\xrightarrow{\mathrm{id}}X$.  Farther out, their descendants may look
unrelated.  One can then ask whether the original identity bridge lifts.  The
lift need not be an isomorphism.  It may be a correspondence, duality,
functor, common invariant, spectral sequence, or categorical equivalence.
Such a relation may initially remain conjectural. Common ancestry guides
the search; it does not by itself prove that the lift exists.

Here \emph{fork} refers to branching developments. 
We reserve \emph{fork object} for a case in which the common base
and lower bridge can be identified precisely.  A \emph{fork candidate} is a
case in which such a structure is suggested but has not yet been made precise.
An initial resemblance is a \emph{forky signal} whose mathematical
content still needs to be examined.

Three-dimensional Floer theory contains several examples.  Flat connections
and representation spaces feed both gauge-theoretic and
symplectic--Lagrangian Floer theories.  Contact three-manifolds lead toward
Reeb dynamics and holomorphic curves on one branch and Seiberg--Witten theory
on another.  Taubes's correspondence identifies embedded contact homology
with an appropriate Seiberg--Witten Floer cohomology \cite{Taubes2010}.  It
shows that descendants which have traveled through very different
mathematical cultures can remain joined by a precise bridge.  This suggests a way to search: determine what each branch preserved
or transformed, record the bridges already known, and look for a structure
on one branch whose counterpart on the other is missing.

A reader who finds this useful can add the following optional entry to the
starter profile.

\begin{profilebox}
\textbf{Optional addition: finding fork objects}

\smallskip
While we work, remain alert to two mathematical developments that may descend
from an overlapping base.  When such a situation appears, identify the shared
objects and relations, the two branches of development, and the bridge visible
at the base---in the simplest case, the identity on a common object.  Then ask
whether this bridge can be lifted progressively through the two theories.
Do not demand an immediate equivalence at the highest level.  Move upward one
step at a time and test whether the bridge survives as a map, correspondence,
functor, duality, common invariant, or other precise relation.  Record where
the lift succeeds, changes form, or meets an obstruction, and distinguish
established bridges from conjectural ones.

Use \emph{fork object} only when the common base and
lower bridge are precise.  Use \emph{fork candidate} when they are plausible
but incomplete, and \emph{forky signal} for a resemblance still to be
examined.  Discard a signal if examination shows that the resemblance has no
useful mathematical content.

If a candidate or forky signal is unrelated to our current problem, do not
interrupt or redirect the work.  Make a short deferred note containing only
the possible common base, the two branches, the bridge at the base, and the
reason a lift might be interesting.  Show me these deferred notices at the end
of the working session.  Preserve a candidate only if I judge it worth
keeping.

With my approval, promising candidates may later be studied by one or more
dedicated agents outside the main session.  Their reports should identify
known bridges, possible lifts, obstructions, and the status of every proposed
connection.  Keep the reports, including failures and their causes, for later
projects. Let me decide what to investigate and preserve.
\end{profilebox}

The instruction can be activated when a possible common base appears.
A short note preserves the candidate until there is time to study it.
If it survives examination, the profile can acquire questions specific
to that fork:
\[
 \begin{aligned}
 \text{core profile}&\longrightarrow\text{fork detected}
 \longrightarrow\text{questions about the fork}\\
 &\longrightarrow\text{focused exploration}.
 \end{aligned}
\]
The underlying move is to descend from developed theories toward a simpler
shared structure, then ask how far the relation can be lifted again.

I tried the instruction while writing this article. Asked for a fork in
symplectic geometry, the model, with my profile available, proposed the
cotangent bundle, leading on one side toward Fukaya categories and on the
other toward constructible sheaves. This is already rich mathematical
territory: the Nadler--Zaslow correspondence \cite{NadlerZaslow2009},
Tamarkin's approach to non-displaceability \cite{Tamarkin2008}, the
Guillermou--Kashiwara--Schapira quantization of Hamiltonian isotopies
\cite{GKS2012}, and the Ganatra--Pardon--Shende work on wrapped Fukaya
categories \cite{GPS2024}.

Near the Riemann hypothesis, it suggested the zeta function, the Gaussian
and theta function, varieties over finite fields, and hyperbolic surfaces.
The suggestions point toward the explicit formula \cite{Riemann1859},
Tate's harmonic analysis \cite{Tate1967}, Weil's conjectures
\cite{Weil1949}, Grothendieck's Lefschetz framework
\cite{Grothendieck1977}, Deligne's proof \cite{Deligne1974}, and Selberg's
trace formula \cite{Selberg1956}. These names locate regions to examine;
they do not yet identify precise forks. We would have to specify the
common base and the first meaningful lift in each case. The prompts
recovered known mathematical connections. Whether a profile helped the
model do so, or could help it find something new, requires a comparison.

Here is an exercise for the reader.  Apply the fork instruction to the work of
Wang and Zahl on the Kakeya conjecture in three dimensions
\cite{WangZahl2022,WangZahl2025}.  What is the common base, which branches grew
from it, which bridges are already known, and which bridge remains to be
lifted?  Our speculative answer was this: the fork object is a family of thin
tubes.  One branch leads to geometric measure theory, union volume, and dimension; the
other to wave packets, Kakeya maximal estimates, Fourier restriction, and
oscillatory integrals.  At the base, the same thin tube is viewed geometrically
and as the physical-space concentration region of a wave packet.  Wang and
Zahl carry the geometric analysis to full dimension through convex
non-clustering, stickiness, grains, and multiscale structure. Their geometric
theorem already permits measurable shadings of the tubes. A further question
is how much of its mechanism can be used in wave-packet estimates, where
amplitudes, phases, orthogonality, and cancellation also matter.

Asked for forks near instanton Floer theory for three-manifolds,
the model returned the Atiyah--Floer
picture.  A Heegaard surface and its space of flat connections form the common
base; one branch develops through gauge theory and instanton Floer homology,
the other through symplectic geometry and Lagrangian Floer homology.  The
conjecture of Atiyah and Floer seeks such a lifted bridge
\cite{Atiyah1988,DaemiFukaya2018}.  Daemi, Fukaya, and Lipyanskiy established an admissible $SO(3)$ version
using suitable splittings with embedded Lagrangians
\cite{DaemiFukayaLipyanskiy2021}. Here the proposed lift already has a
precise mathematical formulation.

For the last trial we supplied some history.  In \emph{The Floer Jungle}, we
describe the Floer learning seminar convened by Jacob Lurie and Akshay
Venkatesh at IAS in 2021.  Venkatesh came away asking whether Floer-type
structures in three- and four-dimensional topology might have ``shadows in
number theory'' \cite{HoferRobertsFloerJungle}.  We then asked:

\begin{quote}
\emph{Are there interesting forks or fork-type structures involving number
theory and near Floer theory?}
\end{quote}

The answer proposed arithmetic representation and deformation spaces as a
possible lower region.  In the topological branch, representations of the
fundamental group appear as flat connections and feed gauge theory, the
Chern--Simons functional, and instanton Floer homology.  In the arithmetic
branch, Galois representations lead to deformation rings and interact with
Selmer groups, Hecke actions, and local--global conditions. The duality
described by Mazur gives one mathematical basis for the analogy between
number fields and three-manifolds
\cite{Mazur1973}; arithmetic gauge and Chern--Simons theories develop parts of
the parallel \cite{Kim2018,KimChernSimons2020}.  Galatius and Venkatesh construct a derived version of Mazur's Galois
deformation ring \cite{GalatiusVenkatesh2020}. The relative Langlands program
of Ben-Zvi, Sakellaridis, and Venkatesh brings in Hamiltonian spaces,
quantization, and conjectural dualities \cite{BZSV2024}.

There is a real limitation in calling this a fork: the
topological and arithmetic representation spaces are analogous rather than a
literally shared base. These works do not provide a full arithmetic
counterpart of the Floer package, including its trajectory-counting
differential, compactness mechanism, and invariance theory. This is better
treated as a fork candidate.  A dedicated agent could test whether a precise
common base and lower bridge can be extracted, whether meaningful partial
lifts exist, or whether the analogy should be dismissed.  Studying such
imperfect cases may also refine the profile entry itself.

The more interesting test is to ask for less familiar bases and incomplete
lifts. A precise obstruction could be as useful as a successful lift if it
identifies where new mathematics is needed. A separate agent could study
such candidates while the main collaboration continues its work. For
example:

\begin{profilebox}
\textbf{Example: an agent searching for fork candidates near RH}

\smallskip
Work outside the main research session.  Search near the Riemann hypothesis
for promising fork candidates other than the zeta function's prime--zero fork,
the Gaussian--theta--Poisson structure, varieties over finite fields, and
hyperbolic surfaces with the Selberg bridge.  Return a short list with reasons for choosing each candidate.

For each candidate, identify the exact common base, the two descendant
branches, the bridge visible at the base, and any known partial lifts, naming
the relevant mathematicians or theorems.  Then formulate the first missing
bridge worth testing and the strongest reason the candidate may be only an
analogy rather than a genuine fork.  Prefer cases in which the higher bridge
region is still sparse enough to offer a possible starting point for a new
project.

Do not interrupt or modify the active research project, and do not present a
candidate as a discovery.  Produce a private candidate report for human review.
Let the researcher decide which entries to keep and which candidates
to investigate further.
\end{profilebox}

History of mathematics acquires an additional importance here. We can
study a discovery to recover questions and ways of thinking that may help
with another problem. A profile can preserve them together with the
circumstances in which they arose. This offers another way for historical scholarship to contribute to
mathematical research. Whether particular entries help is something we can
test.

At IAS, Myles Jackson, Akshay Venkatesh, and I co-direct the History of
Modern Mathematics project \cite{IASHistory2026}. Its interests include
how mathematical ideas are reconceived and rewritten, changing forms of
collaboration, and oral and unpublished mathematical cultures. These are
also places to look for aspects of mathematical thinking that a finished
proof may no longer reveal.

We could also reconstruct a sequence of profiles from a detailed history:
\[
 \Profile_0\longrightarrow\Profile_1\longrightarrow\cdots
 \longrightarrow\Profile_n,
\]
using only what was known at each stage. Testing the profiles on later
problems could show which changes matter. Reproducing a known theorem is
a weak test when the model may have encountered it in training. New
problems or suitably changed examples would tell us more about whether
the profile affects the work.

The history in \emph{The Floer Jungle} suggests how such a reconstruction
might begin \cite{HoferRobertsFloerJungle}. Return to the unresolved
situation before the finished theory made the route seem natural. Which
negative judgments depended on assumptions that could be changed? Which
limited result opened a small door? In the Floer story, the use of flows
could be reconsidered: some essential information survived after machinery
inherited from the surrounding research culture was replaced. Ideas from
distant areas then had to be assembled, with considerable technical work
still ahead.

A profile can retain questions drawn from such a history. This is what we
mean by historical prompting. We want to try the questions on a new problem
and allow the answers to be different. The point is to learn from the
discovery rather than imitate the discoverer.

\begin{profilebox}
\textbf{Optional addition: historical prompting in a mathematical jungle}

\smallskip
When I supply or approve a detailed history of a difficult discovery, recover
from it a research profile that can be tested on a present problem.  Reconstruct
the historical situation before hindsight made the outcome appear natural:
the partial results, failed approaches, competing research cultures, available
tools, authoritative objections, and unresolved technical difficulties.  Do
not imitate the discoverer's personality or retell the finished theory.

When the present problem is messy, look for a small door: a special case,
estimate, failed construction, or limited theorem suggesting that the accepted
picture misses additional structure.  Audit negative judgments at the level of
their exact hypotheses and quantifiers.  Distinguish essential mathematical
information from machinery inherited through research culture, and consider
whether the information can survive after that machinery is replaced.  Search
other fields for tools addressing a named obstruction, stating precisely what
each imported idea would have to accomplish.

Distinguish the proposed idea from the work needed to make it rigorous.
Identify the difficulties in compactness, regularity, transversality,
convergence, invariance, or gluing, and test the most serious one early.  Treat objections, counterexamples, and
failed attempts as information about the problem.  Record provenance and
distinguish what is proved, plausible, and merely hoped for.

When it helps the work, explain what has changed, which opening is
worth pursuing, and what we should test next. Keep the actual mathematical
choices visible rather than turning the history into a heroic story.

\textbf{Kill switch.} If I say ``kill historical prompting,'' ``stop
historical mode,'' or an evident equivalent, stop this mode immediately and
return to the ordinary profile and the immediate task. Reactivate it
only at my request. Recommend stopping if the historical analogy distracts
from the mathematics or supplies a story without a useful test.
\end{profilebox}

Both instructions look beneath a finished theory. Forky prompting returns
to a common base and asks whether a bridge can be lifted. Historical
prompting returns to the circumstances of discovery and asks which
questions can help us in a different situation:
\[
\begin{aligned}
\text{forky:}\quad
 &\text{branches}\longrightarrow\text{base}
 \longrightarrow\text{lifted bridge},\\[3pt]
\text{historical:}\quad
 &\text{finished theory}\longrightarrow\text{live history}
 \longrightarrow\text{new profile}.
\end{aligned}
\]
\section{Learning to work with AI}

A developed profile can bring a mathematician's accumulated understanding
into work with a model. It can also carry a mistake from one problem to
the next. The mathematician must be able to question it, change it, and
occasionally work outside it. This is part of learning to use the new
instrument.

For graduate education, the central task remains to acquire mathematics
deeply enough to ask one's own questions and take responsibility for the
answers. A student should be able to explain an argument, its hypotheses,
and the reason its main idea works. A result obtained with AI can be a
substantial achievement when the student's understanding and choices were
needed to obtain and establish it. An attractive answer to a prompt is
only a beginning if the student cannot yet do that work.

The examples here suggest places to begin. A group checking a paper may
find a mechanism worth developing. A fork may have an incomplete lift or
a structure on one branch with no known counterpart on the other. Such
questions might offer graduate-sized projects, although locating one is
only the first step. The student must learn the surrounding mathematics
and determine whether the proposed problem has substance. A profile can
grow through that process and help carry what was learned into the next
project.

Studying how profiles mark the landscape, interact, and change may lead
to mathematical questions we cannot yet formulate. These need time and
serious examples. Mathematicians can begin this work now, individually,
in research groups, and in seminars, and help shape mathematical work
with AI.

The whole system is dynamic. Most researchers would value the peace of
mind of knowing that they can develop their ideas privately and decide
what to disclose.

\begin{ack}[Use of AI]
After a disappointing trial in 2025, I began sustained mathematical work
with AI in late May 2026. This article grew from that work and was developed
in discussion with OpenAI Codex,
using my research profile and drawing on \cite{Hofer2026}. As of
22 October 2026, the provisional core profile contained 6,610 words,
the collaboration manual 1,509, and two project-specific profile files a
further 4,238, counting the complete files. Detailed mathematical project
records are kept separately.
Codex contributed searches, comparisons, suggestions, drafting, revision,
and typesetting. I am responsible for the article and its mathematical
judgments.
\end{ack}

{\linespread{0.86}
}

\begin{authorinfo}
Helmut H. W. Hofer, professor emeritus at the Institute for Advanced Study,
works in symplectic geometry and mathematical AI.
\email{hofer@math.ias.edu}
\end{authorinfo}

\end{document}